\documentclass[a4paper]{article}

\usepackage{amsmath}
\usepackage{mathtools}
\usepackage{amsfonts}
\usepackage{amssymb}
\usepackage[compress]{cite}
\usepackage{graphicx}
\usepackage{algorithmicx}
\usepackage{algcompatible}
\usepackage{algorithm}
\usepackage{url}
\usepackage{color}

\newtheorem{lem}{\sc Lemma}[section]
\newtheorem{thm}[lem]{\sc Theorem}
\newtheorem{cor}[lem]{\sc Corollary}
\newtheorem{assum}{\bf Assumption}[section]
\newtheorem{rmk}{\bf Remark}[section]
\newtheorem{defn}{\sc Definition}[section]

\newcommand{\QED}{\hfill~\rule[-1pt]{5pt}{5pt}\par\medskip}
\newenvironment{pf}{{\it Proof:\ }}{}

\renewcommand{\matrix}[2]{\left[\begin{array}{#1} #2 \end{array}\right] }

\newcommand{\jpfig}[4]{\begin{figure}[t] \centering \includegraphics[width=#1\linewidth]{#2}  \caption{\label{#3}#4} \end{figure}}
\newcommand{\C}{\mathcal{C}}
\newcommand{\p}{\mathcal{P}}

\newcommand{\K}{\mathcal{K}}

\newcommand{\pp}{\mathbf{X}}

\DeclareMathOperator*{\argmin}{arg\;min}

\DeclareMathOperator*{\trace}{tr}
\DeclareMathOperator*{\diag}{diag}
\DeclareMathOperator*{\esssup}{ess\;sup}

\begin{document}

\title{Adaptive Control Design under Structured Model Information Limitation: A Cost-Biased Maximum-Likelihood Approach\thanks{The work was supported by the Swedish Research Council and the Knut and Alice Wallenberg Foundation.}}

\author{Farhad Farokhi and Karl H. Johansson \thanks{ACCESS Linnaeus Center, School of Electrical Engineering, KTH Royal Institute of Technology, Stockholm, Sweden. Emails:\{farokhi,kallej\}@ee.kth.se}}

\date{}

\maketitle

\begin{abstract} Networked control strategies based on limited information about the plant model usually results in worse closed-loop performance than optimal centralized control with full plant model information. Recently, this fact has been established by utilizing the concept of competitive ratio, which is defined as the worst case ratio of the cost of a control design with limited model information to the cost of the optimal control design with full model information. We show that an adaptive controller, inspired by a controller proposed by Campi and Kumar, with limited plant model information, asymptotically achieves the closed-loop performance of the optimal centralized controller with full model information for almost any plant. Therefore, there exists, at least, one adaptive control design strategy with limited plant model information that can achieve a competitive ratio equal to one. The plant model considered in the paper belongs to a compact set of stochastic linear time-invariant systems and the closed loop performance measure is the ergodic mean of a quadratic function of the state and control input. We illustrate the applicability of the results numerically on a vehicle platooning problem.
\end{abstract}

\section{Introduction}
Networked control systems are often complex large-scale engineered systems, such as power grids~\cite{Massoud1507024}, smart infrastructures~\cite{negenborn2010intelligent}, intelligent transportation systems~\cite{Swaroop917908,Collier272224,Varaiya250509}, or future aerospace systems~\cite{Giulietti887447,Fowler1234649}. These systems consists of several subsystems each one often having many unknown parameters. It is costly, or even unrealistic, to accurately identify all these plant model parameters offline. This fact motivates us to focus on optimal control design under structured parameter uncertainty and limited plant model information constraints.

There are some recent studies in optimal control design with limited plant model information~\cite{Langbort2010,FarokhiPapaerSubmitted1,Farokhi6120356,FarokhiPapaerSubmitted2,Farokhi-thesis2014}. The problem was initially addressed in~\cite{Langbort2010} for designing static centralized controllers for a class of discrete-time linear time-invariant systems composed of scalar subsystem, where control strategies with various degrees of model information were compared using competitive ratio; i.e., the worst case ratio of the cost of a control design with limited model information scaled by the cost of the optimal control design with full model information. The result was generalized to static decentralized controller for a class of systems composed of fully-actuated subsystems of arbitrary order in~\cite{FarokhiPapaerSubmitted1}. More recently, the problem of designing optimal $H_2$ dynamic controllers using limited plant model information was considered in~\cite{Farokhi6120356}. It was shown that, when relying on local model information, the smallest competitive ratio achievable for any control design strategy for distributed linear time-invariant controllers is strictly greater than one; specifically, equal to the square root of two when the $B$-matrix was assumed to be the identity matrix.

In this paper, we generalize the set of applicable controllers to include adaptive controllers. We use the ergodic mean of a quadratic function of the state and control as a performance measure of the closed-loop system. Choosing this closed-loop performance measure allows us to use certain adaptive algorithms available in the literature~\cite{campi:1890,prandini:1499,graves:715,kumar:163}. In particular, we consider an adaptive controller proposed by Campi and Kumar~\cite{campi:1890}\, which uses a cost-biased (i.e., regularized) maximum-likelihood estimator for learning the unknown parts of the model matrices. We prove that this adaptive control design achieves a competitive ratio equal to one and, hence, the smallest competitive ratio that a control design strategy using adaptive controllers can achieve is equal one (since this ratio is always lower-bounded by one). This is contrary to control design strategies that construct linear time-invariant control laws~\cite{Langbort2010,FarokhiPapaerSubmitted1,Farokhi6120356,FarokhiPapaerSubmitted2,Farokhi-thesis2014}. This shows that, although the design of each subcontroller is only relying on local model information, the closed-loop performance can still be as good as the optimal control design strategy with full model information (in the limit).

The rest of the paper is organized as follows. In Section~\ref{sec:PF}, we present the mathematical problem formulation. In Section~\ref{sec:results}, we introduce the Campi--Kumar adaptive controller using only local model information and we show that it achieves a competitive ratio equal to one. We use this adaptive algorithm on a vehicle platooning problem to demonstrate its performance numerically in Section~\ref{sec:Example} and we conclude the paper in Section~\ref{sec:conclusion}.

\subsection{Notation} \label{subsec:notation}
The sets of natural and real numbers are denoted by $\mathbb{N}$ and $\mathbb{R}$, respectively. We define $\mathbb{N}_0=\mathbb{N}\cup \{0\}$. Additionally, all other sets are denoted by calligraphic letters such as $\mathcal{P}$ and $\mathcal{A}$. 

Matrices are denoted by capital roman letters such as $A$. The entry in the $i^{\textrm{th}}$ row and the $j^{\textrm{th}}$ column of matrix $A$ is~$a_{ij}$. Moreover, $A_{ij}$ denotes a submatrix of matrix $A$, the dimension and the position of which will be defined in the text.

$A > (\geq) 0$ means symmetric matrix $A\in \mathbb{R}^{n\times n}$ is positive definite (positive semidefinite) and $A > (\geq) B$ means $A-B > (\geq) 0$. Let $\mathcal{S}_{++}^n$ ($\mathcal{S}_{+}^n$) be the set of symmetric positive definite (positive semidefinite) matrices in $\mathbb{R}^{n\times n}$.

Let matrices $A\in\mathbb{R}^{n\times n}$, $B\in\mathbb{R}^{n\times m}$, $Q\in\mathcal{S}^n_+$, and $R\in\mathcal{S}^m_{++}$ be given such that the pair $(A,B)$ is stabilizable and the pair $(A,Q^{1/2})$ is detectable. We define $\pp(A,B,Q,R)$ as the unique positive definite solution of the discrete algebraic Riccati equation
$$
%P=A^\top P A-A^\top P B \left(B^\top P B+R \right)^{-1}B^\top P A+Q.
X=A^\top X A-A^\top X B \left(B^\top X B+R \right)^{-1}B^\top X A+Q.
$$
In addition, we define
\begin{equation*}
\begin{split}
\mathbf{L}(A,B,Q,R)=-&\left(B^\top \pp(A,B,Q,R) B+R\right)^{-1} B^\top \pp(A,B,Q,R) A.
\end{split}
\end{equation*}
When matrices $Q$ and $R$ are not relevant or can be deduced from the text, we use $\pp(A,B)$ and $\mathbf{L}(A,B)$ instead of $\pp(A,B,Q,R)$ and $\mathbf{L}(A,B,Q,R)$, respectively.

A measurable function $f:\mathcal{Z}\rightarrow \mathbb{R}$ is said to be essentially bounded if there exists a constant $c\in\mathbb{R}$ such that $|f(z)|\leq c$ almost everywhere. The greatest lower bound of these constants is called the essential supremum of $f(z)$, which is denoted by $\esssup_{z\in \mathcal{Z}} f(z)$.

All graphs $G$ considered in this paper are directed with vertex set $\{1,...,N\}$ for a given $N\in\mathbb{N}$. The adjacency matrix $S\in \{0,1\}^{N \times N}$ of $G$ is a matrix whose entry $s_{ij}=1$ if $(j,i) \in E$ and $s_{ij}=0$, otherwise, for all $1\leq i,j\leq N$.

Let mappings $f,g:\mathbb{Z}\rightarrow \mathbb{R}$ be given. Denote $f(k)=O(g(k))$ if $\limsup_{k\rightarrow \infty} |f(k)/g(k)|<\infty$. Similarly, $f(k)=o(g(k))$ if $\limsup_{k\rightarrow \infty} |f(k)/g(k)|=0$.

Finally, $\chi(\cdot)$ denotes the characteristic function, that is, it returns a value equal one if its statement is satisfied and a value equal zero otherwise.

\section{Problem Formulation} \label{sec:PF}
\subsection{Plant Model}
Consider a discrete-time linear time-invariant dynamical system composed of $N$ subsystems, such that the state-space representation of subsystems $i$, $1\leq i\leq N$, is given by
\begin{equation*}
\begin{split}
x_i(k+1)&=\sum_{j=1}^N [A_{ij} x_j(k)+B_{ij} u_j(k)]+w_i(k); \; x_i(0)=0,
\end{split}
\end{equation*}
where $x_i(k)\in\mathbb{R}^{n_i}$, $u_i(k)\in\mathbb{R}^{m_i}$, and $w_i(k)\in\mathbb{R}^{n_i}$ are state, control input, and exogenous input vectors, respectively. We assume that $\{w_i(k)\}_{k=0}^\infty$ are independent and identically distributed Gaussian random variables with zero means $\mathbb{E}\{w_i(k)\}=0$ and unit covariances $\mathbb{E}\{w_i(k)w_i(k)^\top\}=I$. The assumption of unit covariance is without loss of generality and is only introduced to simplify the presentation. To show this, assume that $\mathbb{E}\{w_i(k)w_i(k)^\top\}=H_i\in\mathcal{S}_{++}^{n_i}$ for all $1\leq i\leq N$. Now, using the change of variables $\bar{x}_i(k)=H_i^{-1/2}x_i(k)$ and $\bar{w}_i(k)=H_i^{-1/2}w_i(k)$ for all $1\leq i\leq N$, we get
$$
\bar{x}_i(k+1)\hspace{-.03in}=\hspace{-.03in}\sum_{j=1}^N [\bar{A}_{ij} \bar{x}_j(k)+\bar{B}_{ij} u_j(k)]+\bar{w}_i(k),
$$
in which $\bar{A}_{ij}=H_i^{-1/2}A_{ij}H_j^{1/2}$ and $\bar{B}_{ij}=H_i^{-1/2}B_{ij}$ for all $1\leq i,j\leq N$. This gives $\mathbb{E}\{\bar{w}_i(k)\bar{w}_i(k)^\top\}=I$. In addition, let $w_i(k)$ and $w_j(k)$ be statistically independent for all $1\leq i\neq j\leq N$. Note that this assumption is often justified by the fact that in many large-scale systems, such as smart grids, the subsystems are scattered geographically and, hence, the sources of their disturbances are independent. We introduce the augmented system as
$$
x(k+1)=Ax(k)+Bu(k)+w(k); \; x(0)=0,
$$
where the augmented state, control input, and exogenous input vectors are
\begin{equation*}
\begin{split}
x(k)&=[x_1(k)^\top\;\dots\;x_N(k)^\top]^\top\in\mathbb{R}^n, \\ 
u(k)&=[u_1(k)^\top\;\dots\;u_N(k)^\top]^\top\in\mathbb{R}^m, \\
w(k)&=[w_1(k)^\top\;\dots\;w_N(k)^\top]^\top\in\mathbb{R}^n, \\
\end{split}
\end{equation*}
with $n=\sum_{i=1}^N n_i$ and $m=\sum_{i=1}^N m_i$. In addition, the augmented model matrices are 
$$
B=\matrix{ccc}{B_{11} & \cdots & B_{1N}\\ \vdots & \ddots & \vdots \\ B_{N1} & \cdots & B_{NN}}\in\mathcal{B}\subset\mathbb{R}^{n\times m},
$$ 
and
$$
A=\matrix{ccc}{A_{11} & \cdots & A_{1N}\\ \vdots & \ddots & \vdots \\ A_{N1} & \cdots & A_{NN}}\in\mathcal{A}\subset\mathbb{R}^{n\times n}.
$$
Let a directed plant graph $G_\p$ with its associated adjacency matrix $S^\p$ be given. The plant graph $G_\p$ captures the interconnection structure of the plants, that is, $A_{ij}\neq 0$ only if $s^\p_{ij}\neq 0$. Hence, the sets $\mathcal{A}$ and $\mathcal{B}$ are structured by the plant graph:
\begin{equation*}
\begin{split}
\mathcal{A}\subseteq \bar{\mathcal{A}}=\{A\in\mathbb{R}^{n\times n}&\;|\; s^\p_{ij}= 0 \Rightarrow A_{ij}= 0\in\mathbb{R}^{n_i\times n_j} \mbox{ for all } i,j \mbox{ such that } 1\leq i,j\leq N \},\\
\mathcal{B}\subseteq \bar{\mathcal{B}}=\{B\in\mathbb{R}^{n\times m}&\;|\; s^\p_{ij}= 0 \Rightarrow B_{ij}= 0\in\mathbb{R}^{n_i\times m_j} \mbox{ for all } i,j \mbox{ such that } 1\leq i,j\leq N \},
\end{split}
\end{equation*}
From now on, we present a plant with its pair of corresponding model matrices as $P=(A,B)$ and define $\mathcal{P}=\mathcal{A}\times\mathcal{B}$ as the set of all possible plants. We make the following assumption on the set of all plants:

\begin{assum} \label{asm:1} The set $\mathcal{A}\times \mathcal{B}$ is a compact set (with nonzero Lebesgue measure in the space $\bar{\mathcal{A}}\times \bar{\mathcal{B}}$) and the pair $(A,B)$ is controllable for almost  all $(A,B)\in\mathcal{A}\times\mathcal{B}$. 
\end{assum}

Note that the assumption that the pair $(A,B)$ is controllable for almost all $(A,B)\in\mathcal{A}\times\mathcal{B}$ is guaranteed if and only if the family of systems is structurally controllable~\cite{1100557,Dion20031125}.

\subsection{Adaptive Controller} 
We consider (possibly) infinite-dimensional nonlinear controllers $\mathbf{K}_i=(\mathbf{K}^{(k)}_i)_{k\in\mathbb{N}_0}$ for each subsystem~$i$, $1\leq i\leq N$, with control law
$$
u_i(k)=\mathbf{K}_i^{(k)}(\{x(t)\}_{t=0}^k\cup \{u(t)\}_{t=0}^{k-1}), \hspace{.2in} \forall\;k\in\mathbb{N}_0,
$$
where $\mathbf{K}_i^{(k)}:\prod_{i=1}^k \mathbb{R}^n \times \prod_{i=1}^{k-1} \mathbb{R}^m \rightarrow \mathbb{R}^{m_i} $ is the feedback control law employed at time $k\in\mathbb{N}_0$. Let $\K_i$ denote the set of all such control laws. We also define $\K=\prod_{i=1}^N \K_i$ as the set of all admissible controllers.

\subsection{Control Design Strategy}
A control design strategy $\Gamma$ is a mapping from the set of plants $\mathcal{P}=\mathcal{A}\times \mathcal{B}$ to the set of admissible controllers $\mathcal{K}$. We can partition $\Gamma$ using the control input size as
$$
\Gamma=\matrix{ccc}{\Gamma_{1} \\ \vdots \\ \Gamma_{N}},
$$
where, for each $1\leq i\leq N$, we have $\Gamma_i:\mathcal{A}\times \mathcal{B}\rightarrow \K_i.$ Let a directed design graph $G_\C$ with its associated adjacency matrix $S^\C$ be given. We say that the control design strategy $\Gamma$ satisfies the limited model information constraint enforced by the design graph $G_\C$ if, for all $1\leq i\leq N$, $\Gamma_i$ is only a function of
$$
\{[A_{j1}\;\dots\;A_{jN}],[B_{j1}\;\dots\;B_{jN}]\;|\; s^\C_{ij}\neq 0\}.
$$
The set of all control design strategies that obey the structure given by the design graph $G_\C$ is denoted by $\C$.

\subsection{Performance Metric}
In this paper, we are interested in minimizing the performance criterion
\begin{equation} \label{eqn:cost}
\begin{split}
J_{P}(\mathbf{K})=\limsup_{T\rightarrow \infty} \frac{1}{T} \sum_{k=0}^{T-1} x(k)^\top Q x(k)+u(k)^\top R u(k),
\end{split}
\end{equation}
where $Q\in\mathcal{S}_+^n$ and $R\in\mathcal{S}_{++}^{m}$. We make the following assumption concerning the performance criterion:

\begin{assum} \label{asm:2} The pair $(A,Q^{1/2})$ is observable for almost all $A\in\mathcal{A}$.
\end{assum}

Considering that the observability of the pair $(A,Q^{1/2})$ is equivalent to the controllability of the pair $(A^\top,Q^{1/2})$, we can verify Assumption~\ref{asm:2} using the available results on structural controllability~\cite{1100557,Dion20031125}.

\begin{rmk} Assumptions~\ref{asm:1} and~\ref{asm:2}, that the pair $(A,B)$ is controllable and the pair $(A,Q^{1/2})$ is observable for almost all $(A,B)\in\mathcal{A}\times\mathcal{B}$, originate from the results of Campi and Kumar~\cite{campi:1890}. They used these assumptions to guarantee that the underlying algebraic Riccati equation admits a unique positive-definite solution for almost any selection of model matrices $(A,B)\in\mathcal{A}\times\mathcal{B}$~\cite[p.\,1892]{campi:1890}. We can relax these assumptions for the results in this paper to that the pair $(A,B)$ is stabilizable and the pair $(A,Q^{1/2})$ is detectable for almost all $(A,B)\in\mathcal{A}\times\mathcal{B}$~\cite{1102434}.
\end{rmk}

Note that for linear controllers the performance measure~(\ref{eqn:cost}) represents the $H_2$-norm of the closed-loop system from exogenous input $w(k)$ to output
$$
y(k)=\matrix{c}{Q^{1/2}x(k) \\ R^{1/2}u(k)}.
$$

\begin{defn} Let a plant graph $G_{\mathcal{P}}$ and a design graph $G_{\mathcal{C}}$ be given. Assume that, for every plant $P \in \mathcal{P}$, there exists an optimal controller $\mathbf{K}^*(P) \in \mathcal{K}$ such that
\begin{equation*}
J_P (\mathbf{K}^*(P))\leq J_P (\mathbf{K}), \hspace{.3in} \forall \;\mathbf{K} \in \K.
\end{equation*}
The \emph{average competitive ratio} of a control design method $\Gamma\in\C$ is defined as 
\begin{equation} \label{eqn:acrdef}
r_{\p}^{\mathrm{ave}} (\Gamma)= \int_{\xi\in\p} \frac{J_\xi (\Gamma(\xi))}{J_\xi (\mathbf{K}^*(\xi))} f(\xi)\; \mathrm{d}\xi,
\end{equation}
where $f:\p\rightarrow \mathbb{R}$ is a positive continuous function which shows the relative importance of plants in $\p$. Without loss of generality, we assume that $\int_{\p} f(P) \mathrm{d}P=1$ (up to rescaling $f$ by a constant factor since $\p$ is a compact set and $f$ is a continuous mapping). The \emph{supremum competitive ratio} of a control design method $\Gamma\in\C$ is defined as
\begin{equation} \label{eqn:scrdef}
r_{\p}^{\mathrm{sup}} (\Gamma)= \esssup_{P\in\p} \frac{J_P (\Gamma(P))}{J_P (\mathbf{K}^*(P))}.
\end{equation}
\label{def_comp_rat}
\end{defn}

The mapping $\mathbf{K}^*$ is not required to lie in the set $\C$, and is obtained by searching over the set of centralized controllers with access to the full plant model information. Hence, $\mathbf{K}^*(P)=\mathbf{L}(A,B)$ for all plants $P=(A,B)\in\p$.

The supremum competitive ratio $r_{\p}^{\mathrm{sup}}$ is a modified version of the competitive ratio considered in~\cite{Langbort2010,FarokhiPapaerSubmitted1,Farokhi6120356,FarokhiPapaerSubmitted2,Farokhi-thesis2014}.  Note that using essential supremum in~(\ref{eqn:scrdef}), we are neglecting a subset of plants with zero Lebesgue measure. However, this is not crucial for practical purposes since it is unlikely to encounter such plants in a real situation. As a starting point, let us prove an interesting property relating the average and supremum competitive ratios.

\begin{lem} \label{lem:aveleqsup} For any control design strategy $\Gamma\in\C$, we have
$1 \leq r_\p^{\mathrm{ave}}(\Gamma) \leq r_\p^{\mathrm{sup}}(\Gamma).$
\end{lem}

\begin{pf} See Appendix~\ref{proof:lem:aveleqsup}. \QED\end{pf}

In this paper, we are interested in solving the optimization problem
\begin{equation}
\argmin_{\Gamma\in\C} r_\p(\Gamma), \label{eqn:problem}
\end{equation}
where $r_\p$ is either $r_\p^{\mathrm{ave}}$ or $r_\p^{\mathrm{sup}}$. This problem was studied in~\cite{Farokhi6120356} when the set of plants is fully-actuated discrete-time linear time-invariant systems and the set of admissible controllers is finite-dimensional discrete-time linear dynamic time-invariant systems. It was shown that a modified deadbeat control strategy (which constructs static controllers) is a minimizer of the competitive ratio. Specifically, it was proved that the smallest competitive ratio that a control design strategy which gives decentralized linear time-invariant controllers can achieve is strictly greater than one when relying on local model information. Note that since the optimal control design with full model information is unique (due to Assumption~\ref{asm:2}), even when considering a compact set of plants, the competitive ratio is strictly larger than one for limited model information control design strategies. In this paper, we generalize the formulation of~\cite{Farokhi6120356} to include adaptive controllers. We prove in the next section that we can achieve a competitive ratio equal to one for adaptive controllers. Therefore, we can achieve the optimal performance asymptotically, even if the complete model of the system is not known in advance when designing the subcontrollers.

\section{Main Results} \label{sec:results}
We introduce a specific control design strategy $\Gamma^*$, and subsequently, prove that $\Gamma^*$ is a minimizer of both the average and supremum competitive ratios~$r_\p^{\mathrm{ave}}$ and~$r_\p^{\mathrm{sup}}$. For each plant $P\in\p$, this control design strategy constructs an adaptive controller $\Gamma^*(P)$ using a modified version of the Campi--Kumar adaptive algorithm~\cite{campi:1890}, see Algorithm~\ref{alg:1}. Note that in the Campi--Kumar adaptive algorithm, a central controller estimates the model of the system and controls the system. However, in our modified Campi--Kumar adaptive algorithm in Algorithm~\ref{alg:1}, each subcontroller estimates the model of the system independently and controls its corresponding subsystem separately. Hence, each adaptive subcontroller arrives at different model estimates.

At even time steps in Algorithm~\ref{alg:1}, each subcontroller solves a cost-biased (i.e., regularized) maximum-likelihood problem to extract estimates of the parts of the model matrices that it does not know. In this optimization problem, subcontroller~$i$ fixes the known parts of the model matrices, i.e., $\{[A_{j1}\;\dots\;A_{jN}],[B_{j1}\;\dots\;B_{jN}]|s^\C_{ij}\neq 0\}$, and searches over the unknown parts (see the constraints  in Line~6 of Algorithm~\ref{alg:1}). Due to this information asymmetry, subcontrollers arrive at different model estimates. Upon extracting these estimates, subcontroller~$i$ calculates the optimal control law (by solving the associated Riccati equation) and implements the part that is related to its actuators (see Lines~10 and~11 in Algorithm~\ref{alg:1}).

\begin{rmk} Most often, in practice, some of the entries of the unknown parts of the model matrices are determined by the physical nature of the problem while the rest can vary (due to the parameter uncertainties and the lack of model information from other subsystems). For instance, in heavy-duty vehicle platooning (see Section~\ref{sec:Example}), since the position can ideally be calculated by integrating the velocity over time, some of the entries in the model matrices are fixed (to zero or one). However, other entries may depend on the parameters of the vehicle (e.g., vehicle mass, viscous drag coefficient, and power conversion quality coefficient). Considering that these entries are universally-known constants, one can add them as constraints to the cost-biased maximum-likelihood optimization problem in Algorithm~\ref{alg:1} to reduce the number of decision variables.
\end{rmk}

In Algorithm~\ref{alg:1}, we use the notation $(A^{(i)}(k),B^{(i)}(k))$, at each time step $k\in\mathbb{N}_0$, to denote subsystem~$i$'s estimate of the global system model $P=(A,B)$. Furthermore, for each $1\leq i\leq N$, we use the mapping $\mathbf{T}_i:\mathbb{R}^{m\times n}\rightarrow \mathbb{R}^{m_i\times n}$ defined as
$$
\mathbf{T}_i\matrix{ccc}{X_{11} & \cdots & X_{1N} \\ \vdots & \ddots & \vdots \\ X_{N1} & \cdots & X_{NN} }=\matrix{ccc}{X_{i1} & \cdots & X_{iN} },
$$
where $X_{\ell j}\in\mathbb{R}^{m_\ell\times n_j}$ for each $1\leq \ell,j\leq N$. Let us also, for all $k\in\mathbb{N}_0$, introduce the notation
$$
K(k)=\matrix{c}{\mathbf{T}_1K^{(1)}(k)\\ \vdots \\ \mathbf{T}_NK^{(N)}(k)}\in\mathbb{R}^{m\times n},
$$
where matrices $K^{(i)}(k)$ are defined in Line 10 of Algorithm~\ref{alg:1}. For each $\delta>0$, we introduce
\begin{equation*}
\begin{split}
\mathcal{W}_\delta(A,B):=\{(\bar{A},\bar{B})\in\mathcal{A}\times \mathcal{B} &\;|\; \|[A+B\mathbf{L}(\bar{A},\bar{B})]-[\bar{A}+\bar{B}\mathbf{L}(\bar{A},\bar{B})]\|\geq \delta \}.
\end{split}
\end{equation*}
Let us start by presenting a result on the convergence of the global plant model estimates to the correct value.

\begin{algorithm}[t!]
\caption{Control design strategy $\Gamma^*(P)$.}
\label{alg:1}
\begin{algorithmic}[1]
\STATE \textbf{Parameter:} $\{\mu(k)\}_{k=0}^\infty$ such that $\lim_{k\rightarrow \infty}\mu(k)=\infty$ but $\mu(k)=o(\log(k))$.
\STATE Initialize $(A^{(i)}(0),B^{(i)}(0))$ for all $i\in\{1,\dots,N\}$.
\FOR{$k=1,2,\dots$}
\FOR{$i=1,2,\dots,N$}
\IF{ $k$ is even }
\STATE{Update subsystem~$i$ estimate as
\begin{equation*}
\begin{split}
\hspace{.0in}(A^{(i)}(k),B^{(i)}(k))=&\argmin_{(\hat{A},\hat{B})\in\mathcal{A}\times \mathcal{B}} \mathbf{W}(\hat{A},\hat{B},\mathcal{F}_k), \\ & \mbox{ subject to } \hat{A}_{\ell j}=A_{\ell j}, \hat{B}_{\ell j}=B_{\ell j}, \forall\;j,\ell\in\{1,\dots,N\}, s^\C_{\ell i}\neq 0, \\ & \hspace{.69in} \hat{A}_{z q}=0, \forall\;z,q\in\{1,\dots,N\}, s^\p_{zq}= 0,
\end{split}
\end{equation*}
\hspace{.4in} where
\begin{equation*}
\begin{split}
\mathbf{W}(\hat{A},&\hat{B},\mathcal{F}_k)=\mu(k) \trace (\pp(\hat{A},\hat{B}))+\sum_{t=1}^k\|x(t)-\hat{A}x(t-1)-\hat{B}u(t-1)\|_2^2.
\end{split}
\end{equation*}}
\ELSE
\STATE $(A^{(i)}(k),B^{(i)}(k))\leftarrow (A^{(i)}(k-1),B^{(i)}(k-1))$.
\ENDIF
\STATE $K^{(i)}(k)\leftarrow \mathbf{L}(A^{(i)}(k),B^{(i)}(k))$.
\STATE $u_i(k)\leftarrow\mathbf{T}_iK^{(i)}(k)x(k)$.
\ENDFOR
\ENDFOR
\end{algorithmic}
\end{algorithm}

\begin{lem} \label{tho:1} Let $\Gamma^*(P)$ be defined as in Algorithm~\ref{alg:1} for each plant $P\in\p$. There exists a set $\mathcal{N}\subset\mathcal{P}$ with zero Lebesgue measure (in the space $\bar{\mathcal{A}}\times\bar{\mathcal{B}}$) such that, if $P\notin \mathcal{N}$, then
\begin{equation} \label{eqn:tho:1}
\lim_{k\rightarrow \infty} \pp(A^{(i)}(k),B^{(i)}(k)) \stackrel{as}{\leq} \pp(A,B),
\end{equation}
and
\begin{equation} \label{eqn:tho:2}
\sum_{t=0}^k \chi((A^{(i)}(k),B^{(i)}(k))\in \mathcal{W}_\delta(A,B))\stackrel{as}{=}O(\mu(k)),
\end{equation}
\begin{equation} \label{eqn:tho:2.5}
\sum_{t=0}^k \chi(\|K^{(i)}(k)-\mathbf{L}(A ,B )\|>\rho)\stackrel{as}{=}O(\mu(k)),
\end{equation}
\begin{equation} \label{eqn:tho:3}
\sum_{t=0}^k \chi(\|K(k)-\mathbf{L}(A ,B )\|>\rho)\stackrel{as}{=}O(\mu(k)),
\end{equation}
for all $\delta,\rho>0$, where $x\stackrel{as}{=} y$ and $x\stackrel{as}{\leq} y$ mean $\mathbb{P}\{x= y\}=1$ and $\mathbb{P}\{x \leq y\}=1$, respectively. In addition, we get
\begin{equation} \label{eqn:tho:4}
\limsup_{T\rightarrow\infty} \frac{1}{T} \sum_{k=0}^{T-1} \|x(k)\|^p+\|u(k)\|^p \stackrel{as}{<}\infty,\; \forall\;p\geq 1.
\end{equation}
\end{lem}

\begin{pf} See Appendix~\ref{proof:tho:1}. \QED\end{pf}

Note that, according to Lemma~\ref{tho:1}, we know that there exists a set $\mathcal{N}\subset\mathcal{P}$ with zero Lebesgue measure such that, if $P\notin \mathcal{N}$, the estimates in the modified Campi--Kumar adaptive algorithm (Algorithm~\ref{alg:1}) converge to the correct global plant model. This fact is a direct consequence of using regularized maximum-likelihood estimators in the Campi--Kumar algorithm~\cite{Kumar52293}. We need the following lemma.

\begin{lem} \label{lem:norm} For any matrices $X,P,Y\in\mathbb{R}^{n\times n}$, we have
$$
\|X^\top P X-Y^\top P Y\|\leq \|P\| \|X-Y\|(\|X\|+\|Y\|).
$$
\end{lem}

\begin{pf} See Appendix~\ref{proof:lem:norm}. \QED\end{pf}

Now, we are ready to present the main result of this section.

\begin{thm}\label{tho:3} Let $\Gamma^*(P)$ be defined as in Algorithm~\ref{alg:1} for each plant $P\in\p$. There exists a set $\mathcal{N}\subset\mathcal{P}$ with zero Lebesgue measure such that, if $P\notin \mathcal{N}$, then
$$
J_P(\Gamma^*(P)) \stackrel{as}{=}J_P(\mathbf{K}^*(P)).
$$
\end{thm}

\begin{pf} The proof follows the same reasoning as in~\cite{campi:1890}. The main difference between the proof of this theorem and that of Theorem~6 in~\cite{campi:1890} is caused by that the subcontrollers creates local model estimates $(A^{(i)}(k),B^{(i)}(k))$ and local control gains $K^{(i)}(k)$, which can technically be different from each other (because they rely on their private information). Moreover, $K^{(i)}(k)$ is the control gain that subcontroller~$i$ creates for the entire system (through solving the underlying Riccati equation based on its own model estimates); however, it can only use its corresponding actuators to implement the control signal $\mathbf{T}_iK^{(i)}(k)x(k)$. Therefore, considering that all these local control gains (and their corresponding closed-loop systems) approach each other (see Lemma~\ref{tho:1}) and, ultimately, converge to the true optimal controller, one needs to show that their contribution to the cost function $J$ also converges to that of the optimal controller with full model information (see Items~3 and~5 from the proof below). 

According to~\cite[p.158]{kumar1986stochastic}, for all $1\leq i\leq N$, we get the set of equations in
\begin{align} 
\trace\{\pp(A^{(i)}&(k),B^{(i)}(k))\}+x(k)^\top \pp(A^{(i)}(k),B^{(i)}(k)) x(k) \nonumber
\\&=x(k)^\top Q x(k)+u^{(i)}(k)^\top R u^{(i)}(k) +\mathbb{E}\{
(A^{(i)}(k)x(k)+B^{(i)}(k)u^{(i)}(k)+w(k))^\top \nonumber\\&\hspace{.4in} \times\pp(A^{(i)}(k),B^{(i)}(k)) (A^{(i)}(k)x(k)+B^{(i)}(k)u^{(i)}(k)+w(k))\;|\; \mathcal{F}_{k-1}\} \nonumber\\&=x(k)^\top Q x(k)+u^{(i)}(k)^\top R u^{(i)}(k) +\mathbb{E}\{x(k+1)^\top \pp(A^{(i)}(k),B^{(i)}(k)) x(k+1) \;|\; \mathcal{F}_{k-1} \} \nonumber\\&\hspace{.4in} +(A^{(i)}(k)x(k)+B^{(i)}(k)u^{(i)}(k))^\top \pp(A^{(i)}(k),B^{(i)}(k))(A^{(i)}(k)x(k)+B^{(i)}(k)u^{(i)}(k)) \nonumber\\&\hspace{.4in} - (A x(k)+B u(k))^\top \pp(A^{(i)}(k),B^{(i)}(k)) (A x(k)+B u(k)),\label{eqn:proof:longeq}
\end{align}
with $u^{(i)}(k)=K^{(i)}(k)x(k)$ and $u(k)=K(k)x(k)$. Averaging both sides of~(\ref{eqn:proof:longeq}) over time and all subsystems, we get
\begin{align}
\zeta_1(T)+\frac{1}{NT}\sum_{k=0}^{T-1}&\sum_{i=1}^N x(k)^\top \pp(A^{(i)}(k),B^{(i)}(k)) x(k)
\nonumber\\=&\frac{1}{T}\sum_{k=0}^{T-1}x(k)^\top Q x(k)+\frac{1}{NT}\sum_{k=0}^{T-1}\sum_{i=1}^N u^{(i)}(k)^\top R u^{(i)}(k)\nonumber\\&+\frac{1}{NT}\sum_{k=0}^{T-1}\sum_{i=1}^N \mathbb{E}\{x(k+1)^\top \pp(A^{(i)}(k),B^{(i)}(k)) x(k+1) \;|\; \mathcal{F}_{k-1} \}+\zeta_2(T), \label{eqn:proof:longeq:1}
\end{align}
where
$$
\zeta_1(T)=\frac{1}{NT}\sum_{k=0}^{T-1}\sum_{i=1}^N\trace\{\pp(A^{(i)}(k),B^{(i)}(k))\},
$$
and $\zeta_2(T)$ is given in
\begin{equation} \label{eqn:proof:longeq:4}
\begin{split}
\zeta_2(T)=\frac{1}{NT}\sum_{k=0}^{T-1}\sum_{i=1}^N&(A^{(i)}(k)x(k)+B^{(i)}(k)u^{(i)}(k))^\top \pp(A^{(i)}(k),B^{(i)}(k))(A^{(i)}(k)x(k)+B^{(i)}(k)u^{(i)}(k))\\&\hspace{.7in} - (A x(k)+B u(k))^\top \pp(A^{(i)}(k),B^{(i)}(k)) (A x(k)+B u(k)).
\end{split}
\end{equation}
Moreover, $\mathcal{F}_k=\{x(t)\}_{t=0}^k\cup \{u(t)\}_{t=0}^{k-1}$ denotes the observation history.
Subtracting the term 
$\frac{1}{NT}\sum_{k=0}^{T-1}\sum_{i=1}^N \mathbb{E}\{x(k+1)^\top \pp(A^{(i)}(k+1),B^{(i)}(k+1))  x(k+1) |\mathcal{F}_{k-1} \}$
from both sides of~(\ref{eqn:proof:longeq:1}) while adding and subtracting $\frac{1}{T}\sum_{k=0}^{T-1}  u(k)^\top R u(k)$
from right-hand side of~(\ref{eqn:proof:longeq:1}), we get
\begin{equation} \label{eqn:proof:longeq:2}
\begin{split}
\frac{1}{T}\sum_{k=0}^{T-1}[x(k)^\top Q &x(k)+u(k)^\top R u(k)]+\zeta_4(T)+\zeta_5(T)+\zeta_2(T)=\zeta_1(T)+\zeta_3(T),
\end{split}
\end{equation}
where
\begin{equation*}
\begin{split}
&\zeta_3(T)=\frac{1}{NT}\sum_{k=0}^{T-1}\sum_{i=1}^N x(k)^\top \pp(A^{(i)}(k),B^{(i)}(k))x(k) \\& \hspace{1in}-\mathbb{E}\{x(k+1)^\top \pp(A^{(i)}(k\hspace{-.03in}+1),B^{(i)}(k\hspace{-.03in}+1))x(k+1) \;|\; \mathcal{F}_{k-1} \},
\end{split}
\end{equation*}
$$
\zeta_4(T)=\frac{1}{NT}\sum_{k=0}^{T-1}\sum_{i=1}^N u^{(i)}(k)^\top R u^{(i)}(k)-u(k)^\top R u(k),
$$
and
\begin{equation*}
\begin{split}
\zeta_5(T)=&\frac{1}{NT}\sum_{k=0}^{T-1}\sum_{i=1}^N \mathbb{E}\{x(k+1)^\top [\pp(A^{(i)}(k),B^{(i)}(k))-\pp(A^{(i)}(k+1),B^{(i)}(k+1))] x(k+1) \;|\; \mathcal{F}_{k-1} \}.
\end{split}
\end{equation*}
In the rest of the proof, we study the asymptotic behavior of the sequences $\{\zeta_\ell(k)\}_{k=0}^\infty$ for all $1\leq \ell\leq 5$.
\par \noindent \textsc{Item~1}: Asymptotic behavior of $\zeta_1(T)$: First, note that
\begin{align}
\limsup_{T\rightarrow \infty}\zeta_1(T)&=\limsup_{T\rightarrow \infty}\frac{1}{NT} \sum_{k=0}^{T-1} \sum_{i=1}^N \trace\{\pp(A^{(i)}(k),B^{(i)}(k))\} \nonumber
\\&=\frac{1}{N}\sum_{i=1}^N \limsup_{T\rightarrow \infty}\frac{1}{T} \sum_{k=0}^{T-1} \trace\{\pp(A^{(i)}(k),B^{(i)}(k))\}. \label{eqn:tho:revision:1}
\end{align}
Using~(\ref{eqn:tho:1}) inside~\eqref{eqn:tho:revision:1}, we get
\begin{equation*}
\begin{split}
\limsup_{T\rightarrow \infty}\zeta_1(T) \stackrel{as}{\leq} \trace\{\pp(A ,B )\}.
\end{split}
\end{equation*}
\vspace{-.2in}
\par \noindent \textsc{Item~2}: Asymptotic behavior of $\zeta_3(T)$: With a similar strategy as in case~(B) in the proof of Theorem~6 in~\cite{campi:1890}, we can prove that
\begin{equation*}
\begin{split}
&0\stackrel{as}{=}\limsup_{T\rightarrow \infty}\frac{1}{T}\sum_{k=0}^{T-1} x(k)^\top \pp(A^{(i)}(k),B^{(i)}(k)) x(k)\\&\hspace{1in}-\mathbb{E}\{x(k+1)^\top \pp(A^{(i)}(k+1),B^{(i)}(k+1)) x(k+1) \;|\; \mathcal{F}_{k-1} \}.
\end{split}
\end{equation*}
Hence, $\limsup_{T\rightarrow \infty} \zeta_3(T)\stackrel{as}{=}0$.
\par \noindent \textsc{Item~3}: Asymptotic behavior of $\zeta_4(T)$: In this case, we have
\begin{equation*}
\begin{split}
\left|\frac{1}{T}\sum_{k=0}^{T-1}u^{(i)}(k)^\top R u^{(i)}(k)-u(k)^\top R u(k)\right|  &\leq
\frac{1}{T}\sum_{k=0}^{T-1}\left|u^{(i)}(k)^\top R u^{(i)}(k)-u(k)^\top R u(k)\right|
\\ &\leq \frac{1}{T}\sum_{k=0}^{T-1}\|K^{(i)}(k)^\top R K^{(i)}(k)-K(k)^\top R K(k)\| \|x(k)\|^2.
\end{split}
\end{equation*}
According to Lemma~\ref{lem:norm}, we have $\|K^{(i)}(k)^\top R K^{(i)}(k)-K(k)^\top R K(k)\|\leq \|R\|\|K^{(i)}(k)-K(k)\|\linebreak \times (\|K^{(i)}(k)\|+\|K(k)\|)$. Considering that $\mathbf{L}(\cdot,\cdot)$ is a continuous function of its arguments (see~\cite{Delchamps4047792}) and $\p$ is a compact set, we know that $\|K^{(i)}(k)\|$ and $\|K(k)\|$ are uniformly bounded. Hence, $\|K^{(i)}(k)\|+\|K(k)\| \leq M$. Now, using Cauchy--Schwartz inequality~\cite[p.\,98]{friedman1970foundations}, we get
\begin{equation*}
\begin{split}
\left|\frac{1}{T}\sum_{k=0}^{T-1}u^{(i)}(k)^\top R u^{(i)}(k)-u(k)^\top R u(k)\right|^2 &
\leq\hspace{-.02in}\|R\| M \hspace{-.06in} \left(\frac{1}{T}\sum_{k=0}^{T-1}\|K^{(i)}(k)-K(k)\|^2\right) \hspace{-.08in} \left(\frac{1}{T}\sum_{k=0}^{T-1}\left\|x(k) \right\|^4\right)\hspace{-.03in}.
\end{split}
\end{equation*}
Let us introduce the notation $K^o=\mathbf{L}(A,B)$. Note that, for all $\rho>0$, we have
\begin{equation*}
\begin{split}
\frac{1}{T}\sum_{k=0}^{T-1}\|K^{(i)}(k)-K^o\|^2
\leq \rho^2 \hspace{-.05in}+\frac{1}{T}\sum_{k=0}^{T-1} \|K^{(i)}(k)-K^o\|^2
\chi(\|K^{(i)}(k)-K^o\|>\rho).
\end{split}
\end{equation*}
Again, considering the facts that $\mathbf{L}(\cdot,\cdot)$ is a continuous function of its arguments and $\p$ is a compact set, we know that $\|K^{(i)}(k)-K^o\|$ is uniformly bounded. Thus, using~(\ref{eqn:tho:2.5}) from Lemma~\ref{tho:1}, we can show that $\limsup_{T\rightarrow \infty}\frac{1}{T}\sum_{k=0}^{T-1}\|K^{(i)}(k)-K^o\|\stackrel{as}{\leq}\rho^2,$ for all $\rho>0$. Since the choice of $\rho$ was arbitrary, we get $\limsup_{T\rightarrow \infty}\frac{1}{T}\sum_{k=0}^{T-1} \|K^{(i)}(k)-K^o\|^2\stackrel{as}{=}0.$ With a similar reasoning, we can also prove that $\limsup_{T\rightarrow \infty}\frac{1}{T}\sum_{k=0}^{T-1}\|K(k)-K^o\|^2\stackrel{as}{=}0.$ Therefore, considering that $\|K^{(i)}(k)-K(k)\|^2 \leq \|K^{(i)}(k)-K^o\|^2+ \|K(k)-K^o\|^2$, we have $\limsup_{T\rightarrow \infty}\frac{1}{T} \sum_{k=0}^{T-1} \|K^{(i)}(k)- K(k)\|^2\stackrel{as}{=}0$. Hence, $\limsup_{T\rightarrow \infty} \zeta_4(T)\stackrel{as}{=}0$ due to the fact that $\limsup_{T\rightarrow \infty} \left\|x(k) \right\|^4 \stackrel{as}{<}\infty$ according to~(\ref{eqn:tho:4}).
\par \noindent \textsc{Item~4}: Asymptotic behavior of $\zeta_5(T)$: With the same approach as in case~(C) in the proof of Theorem~6 in~\cite{campi:1890}, we can prove
$\limsup_{T\rightarrow \infty} \zeta_5(T)\stackrel{as}{=}0$.

\par \noindent \textsc{Item~5}: Asymptotic behavior of $\zeta_2(T)$: Let us start with studying the asymptotic behavior of the sequence $\{\hat{\zeta}^{(i)}_2(T)\}_{T=0}^\infty$ in
\begin{align} 
\hat{\zeta}^{(i)}_2(T)=\frac{1}{T}\sum_{k=0}^{T-1}&x(k)^\top(A^{(i)}(k)+B^{(i)}(k)K^{(i)}(k))^\top \pp(A^{(i)}(k),B^{(i)}(k))(A^{(i)}(k)+B^{(i)}(k)K^{(i)}(k))x(k)\nonumber\\&\hspace{.7in} - x(k)^\top(A +B K(k))^\top \pp(A^{(i)}(k),B^{(i)}(k)) (A +B K(k))x(k).\label{eqn:proof:longeq:5}
\end{align}
Using Lemma~\ref{lem:norm}, we can upper bound each term as in
\begin{align} 
x(k)^\top(A^{(i)}(k)+&B^{(i)}(k)K^{(i)}(k))^\top \pp(A^{(i)}(k),B^{(i)}(k))(A^{(i)}(k)+B^{(i)}(k)K^{(i)}(k))x(k) \nonumber\\& - x(k)^\top(A +B K(k))^\top \pp(A^{(i)}(k),B^{(i)}(k)) (A +B K(k))x(k) \nonumber \\  &\hspace{-0.4in}\leq\;
\left\|x(k)\right\|\left\|\pp(A^{(i)}(k),B^{(i)}(k))\right\|
\left\|[A^{(i)}(k)+B^{(i)}(k)K^{(i)}(k)] -[A +B K(k)]\right\|
 \nonumber\\&\hspace{.7in}\times \left\|[A^{(i)}(k)+B^{(i)}(k)K^{(i)}(k)]+[A +B K(k)]\right\|. \label{eqn:proof:longeq:6}
\end{align}
Considering again that $\mathbf{L}(\cdot,\cdot)$ and $\pp(\cdot,\cdot)$ are continuous functions of their arguments (see~\cite{Delchamps4047792}) and $\p$ is a compact set, we know that
\begin{equation*}
\begin{split}
\left\|\pp(A^{(i)}(k),B^{(i)}(k))\right\|&\leq M_1,\\
\left\|[A^{(i)}(k)+B^{(i)}(k)K^{(i)}(k)]+[A +B K(k)]\right\| &\leq M_2.
\end{split}
\end{equation*}
Using Cauchy--Schwartz inequality, we get the inequality in
\begin{equation} \label{eqn:proof:longeq:7}
\begin{split}
\hat{\zeta}^{(i)}_2(T)&\leq M_1M_2\frac{1}{T}\sum_{k=0}^{T-1}\left\|x(k)\right\|^2\left\|[A^{(i)}(k)+B^{(i)}(k)K^{(i)}(k)]
-[A +B K(k)]\right\|\\& \leq M_1M_2 \left(\frac{1}{T}\sum_{k=0}^{T-1}\left\|x(k)\right\|^4\right)^{1/2}
\left( \frac{1}{T}\sum_{k=0}^{T-1}\left\|[A^{(i)}(k)+B^{(i)}(k)K^{(i)}(k)]
-[A +B K(k)]\right\|^2\right)^{1/2}.
\end{split}
\end{equation}
Now, note that
\begin{equation*}
\begin{split}
\left\|[A^{(i)}(k)+B^{(i)}(k)K^{(i)}(k)]-[A +B K(k)]\right\|^2 & \leq
\left\|[A^{(i)}(k)+B^{(i)}(k)K^{(i)}(k)]-[A +B K^{(i)}(k)]\right\|^2
\\&\hspace{.5in}+\|[A +B K^{(i)}(k)]-[A +B K^o ]\|^2
\\&\hspace{.5in}+\|[A +B K^o ]-[A +B K(k)]\|^2
\\  & \hspace{.1in}\leq \left\|[A^{(i)}(k)+B^{(i)}(k)K^{(i)}(k)]-[A +B K^{(i)}(k)]\right\|^2
\\&\hspace{.5in}+\|B \|^2\left(\|K^{(i)}(k)-K^o \|^2+\|K(k)-K^o \|^2 \right).
\end{split}
\end{equation*}
Hence, with similar argument as above, we can prove that $\limsup_{T \rightarrow \infty} \hat{\zeta}^{(i)}_2(T)\stackrel{as}{=}0,$ and as a result $\limsup_{T \rightarrow \infty} \zeta_2(T)=\limsup_{T \rightarrow \infty} \frac{1}{N}\sum_{i=1}^N \hat{\zeta}^{(i)}_2(T)\stackrel{as}{=}0.$
Now, we are ready to prove the statement of this theorem. From the asymptotic behavior of sequences $\zeta_1(T)$ and $\zeta_3(T)$, we know that
\begin{equation} \label{eqn:proof:intermediate:inequality:1}
\begin{split}
\trace\{\pp(A, &B)\}\stackrel{as}{\geq} \limsup_{T \rightarrow \infty} \zeta_1(T)+\zeta_3(T).
\end{split}
\end{equation}
Using identity~(\ref{eqn:proof:longeq:2}) inside inequality~(\ref{eqn:proof:intermediate:inequality:1}) shows that
\begin{equation*}
\begin{split}
\trace\{\pp(A,&B)\}\stackrel{as}{\geq}\limsup_{T \rightarrow \infty} \zeta_4(T)+\zeta_5(T)+\zeta_2(T)+\limsup_{T \rightarrow \infty} \frac{1}{T}\sum_{k=0}^{T-1} [x(k)^\top Q x(k)
+u(k)^\top R u(k)],
\end{split}
\end{equation*}
which result in
\begin{equation*}
\begin{split}
\trace\{\pp(A,B)\} \stackrel{as}{\geq}\limsup_{T \rightarrow \infty} \frac{1}{T}\sum_{k=0}^{T-1} [x(k)^\top Q x(k)+u(k)^\top R u(k)].
\end{split}
\end{equation*}
This inequality finishes the proof.\QED
\end{pf}

Now, we are ready to present the solution of problem~(\ref{eqn:problem}).

\begin{cor} \label{cor:1} For any plant graph $G_\p$ and design graph $G_\C$, we get
$r_\p^{\mathrm{ave}}(\Gamma^*)\stackrel{as}{=} 1$ and $r_\p^{\mathrm{sup}}(\Gamma^*)\stackrel{as}{=}1. $
\end{cor}

\begin{pf} See Appendix~\ref{proof:cor:1}. \QED \end{pf}

Corollary~\ref{cor:1} shows that, irrespective of the plant graph $G_\p$ and design graph $G_\C$, there exists a limited model information control design strategy that can achieve a competitive ratio equal one. This control design strategy gives adaptive controllers achieving asymptotically the closed-loop performance of optimal control design strategy with full model information. Note that earlier results stated that such competitive ratio cannot be achieved by static or linear time-invariant dynamic controllers~\cite{Langbort2010,FarokhiPapaerSubmitted1,Farokhi6120356,FarokhiPapaerSubmitted2,Farokhi-thesis2014}.

\section{Example} \label{sec:Example}
As a simple numerical example, let us consider the problem of regulating the distance between~$N$ vehicles in a platoon. We model vehicle~$i$, $1\leq i\leq N$, as
\begin{equation*}
\begin{split}
\matrix{c}{x_i(k+1) \\ v_i(k+1)}=&\left(I+ \Delta T \matrix{cc}{0 & 1 \\ 0 & -\alpha_i/m_i } \right) \matrix{c}{x_i(k) \\ v_i(k)}+\matrix{c}{0 \\ \Delta T \beta_i/m } \bar{u}_i(k)+\matrix{c}{\bar{w}_1^{i}(k) \\ \bar{w}_2^{i}(k)},
\end{split}
\end{equation*}
where $x_i(k)$ is the vehicle's position, $v_i(k)$ its velocity, $m_i$ the mass, $\alpha_i$ the viscous drag coefficient, $\beta_i$ the power conversion quality coefficient, and $\Delta T$ the sampling time. For each vehicle, stochastic exogenous inputs $\bar{w}^i_j(k)\in\mathbb{R}^n$, $j=1,2$, capture the effect of wind, road quality, friction, etc. A discussion regarding the modeling can be found in~\cite{6426380}. For simplicity of presentation, let us consider the case of $N=2$ vehicles. In addition, assume that $\Delta T=1$. As performance objective, the designer wants to minimize the cost function
\begin{equation*}
\begin{split}
J=\limsup_{T\rightarrow \infty}&\frac{1}{T}\sum_{k=0}^{T-1} \big[q_d(x_1(k)-x_2(k)-d^*)^2 +\sum_{i=1,2} q_v(v_i(k)-v^*)^2 + r(\bar{u}_i(k)-\bar{u}_i^*)^2\big],
\end{split}
\end{equation*}
where $q_d$, $q_v$, and $r$ are positive constants that adjust the penalty terms on the position error, the velocity errors, and the control actions. Moreover, $d^*$ and $v^*$ denote the desired distance and velocity of the platoon. Through minimizing $J$, we can regulate the distance between the trucks and their velocity using the least amount of control effort. Note that $\bar{u}_i^*=\alpha_i v^*/\beta_i$ is the average control signal. We can write the reduced-order system using the distance between vehicles and their velocities as state variables in the form
\begin{equation} \label{eqn:vehicledynamics}
z(k+1)=Az(k)+Bu(k)+w(k),\;z(0)=0,
\end{equation}
where
\begin{equation*}
\begin{split}
&z(k)=[v_1(k)-v^* \;,\; x_1(k)-x_2(k)-d^* \;,\; v_2(k)-v^* ]^\top, \\ & u(k)=[\bar{u}_1(k)-\bar{u}_1^* \;,\; \bar{u}_2(k)- \bar{u}_2^* ]^\top, \\ & w(k)=[\bar{w}_2^1(k)\;,\; \bar{w}_1^1(k)+\bar{w}_1^2(k) \;,\; \bar{w}_2^2(k)]^\top,
\end{split}
\end{equation*}
and
\begin{equation*}
\begin{split}
A=\matrix{cccccc}{ 1-\frac{\alpha_1}{m_1} & 0 & 0 \\ 1 & 1 & -1 \\ 0 & 0 & 1- \frac{\alpha_2}{m_2}}, \hspace{.2in} B=\matrix{cc}{ \frac{\beta_1}{m_1} & 0 \\ 0 & 0 \\ 0 & \frac{\beta_2}{m_2} }.
\end{split}
\end{equation*}
This model leads to
\begin{equation} \label{eqn:costvehicle}
J=\limsup_{T\rightarrow \infty}\frac{1}{T}\sum_{k=0}^{T-1} z(k)^TQz(k)+u(k)^TRu(k),
\end{equation}
where $Q=\diag(q_v,q_d,q_v)$ and $R=\diag(r,r)$. To simplify the presentation, let $Q=I$ and $R=I$.

Note that $z(0)=0$ in~(\ref{eqn:vehicledynamics}) indicates that the vehicles start at the desired distance $d^*$ of each other and with velocity $v^*$. However, due to the exogenous inputs $w(k)$, the vehicles drift away from this ideal situation. By minimizing the closed-loop performance criterion in~(\ref{eqn:costvehicle}), the designer minimizes this drift using the least amount of control effort possible.

We define the first subsystem as $\underline{z}_1(k)=z_1(k)$ and the second subsystem as $\underline{z}_2(k) =[z_2(k) \; z_3(k) ]^T$. Therefore, we get
$$
\underline{z}_1(k+1)=a_{11}\underline{z}_1(k)+b_{11}u_1(k)+w_1(k),
$$
and
\begin{equation*}
\begin{split}
\underline{z}_2(k+1)=&\matrix{c}{1\\0}\underline{z}_1(k)  +
\matrix{cc}{1 & -1\\0 & a_{22}}\underline{z}_2(k) + 
\matrix{c}{0 \\ b_{22}} u_2(k)+\matrix{c}{w_2(k) \\ w_3(k)}\hspace{-.03in},
\end{split}
\end{equation*}
where $(a_{ii},b_{ii})$ are local parameters of subsystem $i$. Assume that
$$
\mathcal{A}=\left\{A\in\mathbb{R}^{3\times 3} \;\bigg|\; A=\matrix{cccccc}{ \hspace{-.06in}a_{11} & 0 & 0\hspace{-.06in} \\ \hspace{-.06in}1 & 1 & -1\hspace{-.06in} \\ \hspace{-.06in}0 & 0 & a_{22}\hspace{-.06in}}, a_{11},a_{22}\in [0,1] \right\},
$$
$$
\mathcal{B}=\left\{B\in\mathbb{R}^{3\times 2}\;\bigg|\; B=\matrix{cc}{\hspace{-.06in}b_{11} & 0\hspace{-.06in} \\ \hspace{-.06in}0 & 0\hspace{-.06in} \\ \hspace{-.06in}0 & b_{22}\hspace{-.06in}},   b_{11},b_{22}\in [0.5,1.5] \right\}.
$$

We compare the performance of the introduced adaptive controller with a deadbeat control design strategy $\Gamma^\Delta:\mathcal{P}\rightarrow \mathbb{R}^{2\times 3}$ for this special family of systems as
$$
\Gamma^\Delta(P)=\matrix{ccc}{-a_{11}/b_{11} & 0 & 0 \\ 1/b_{22} & 1/b_{22} & -(1+a_{22})/b_{22}},
$$
for all $P=(A,B)\in\mathcal{P}$. Note that $\Gamma^\Delta$ is a limited model information control design strategy, because each local controller~$i$ is based on only parameters of subsystem~$i$, $i=1,2$. We also compare the results with the centralized Campi--Kumar adaptive controller $\Gamma^{\mathrm{C}}(P)$ in~\cite{campi:1890}. Notice that this control design strategy does not use the model information that is already available to each local controller.

Figure~\ref{figure1} illustrates the running cost of the closed-system with the optimal control design with full model information $\mathbf{K}^*(P)$ (solid red curve), the modified Campi--Kumar adaptive controller $\Gamma^*(P)$ (dashed green curve), the deadbeat control design strategy $\Gamma^\Delta(P)$ (dotted black curve), and the centralized Campi--Kumar adaptive controller $\Gamma^{\mathrm{C}}(P)$ (dashed-dotted magenta curve). The running costs of the closed-system with the modified Campi--Kumar adaptive controller $\Gamma^*(P)$, the centralized Campi--Kumar adaptive controller $\Gamma^{\mathrm{C}}(P)$, and the optimal control design with full model information $\mathbf{K}^*(P)$ both converge to $\trace\{\pp(A,B)\}$ (the horizontal line) as time goes to infinity. The cost of the optimal control design strategy with global model knowledge is always lower than the cost of the adaptive controllers. Moreover, the cost of the modified Campi--Kumar adaptive controller $\Gamma^*(P)$ is always lower than the centralized Campi--Kumar adaptive controller $\Gamma^{\mathrm{C}}(P)$ because $\Gamma^*(P)$ uses the private model information that is available is each local controller, however, $\Gamma^{\mathrm{C}}(P)$ ignores this information. The simulation is done for randomly-selected parameters $(a_{11},b_{11})=(0.4360,1.0497)$ and $(a_{22},b_{22})=(0.0259,0.9353)$. Figure~\ref{figure2} illustrates the convergence of the individual model parameters $(a_{ii},b_{ii})$, $i=1,2$, for the adaptive subcontrollers. Note that only one of the subsystems needs to estimate each parameter (as each one has access to its own model parameters). Moreover, the results of Lemma~\ref{tho:1} imply that the number of instances that the parameter estimation error is above a fixed threshold grows logarithmically. Therefore, such occurrences become rarer in average. However, this does not imply that at any given time, or even on any finite horizon, the estimation error is decreasing as one may notice from $|b_{22}-b_{22}^{(1)}(k)|$ (the dashed-dotted line) in Figure~\ref{figure2}. 

\jpfig{0.45}{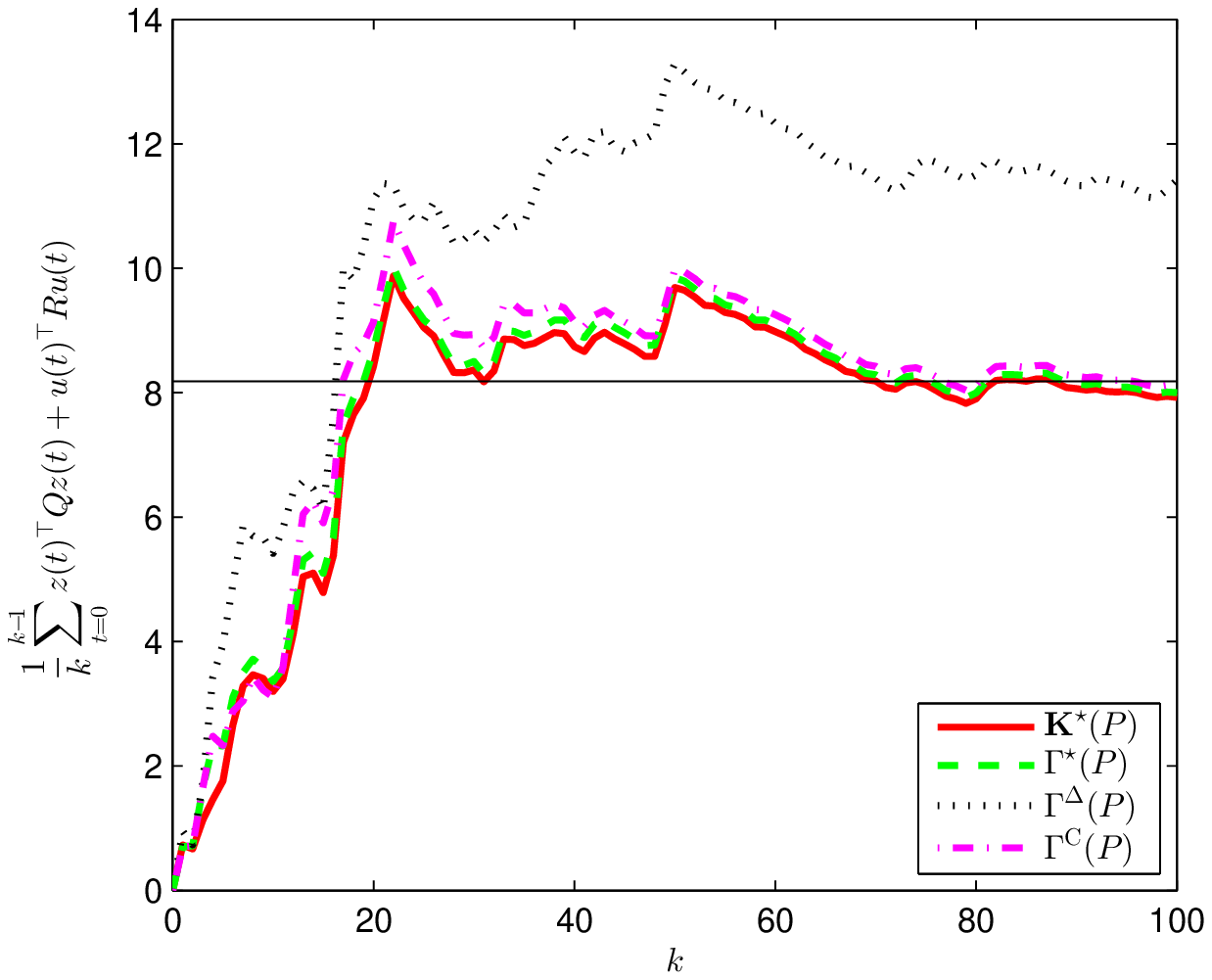}{figure1}{The running cost of the closed-system for four controllers. }

\jpfig{0.45}{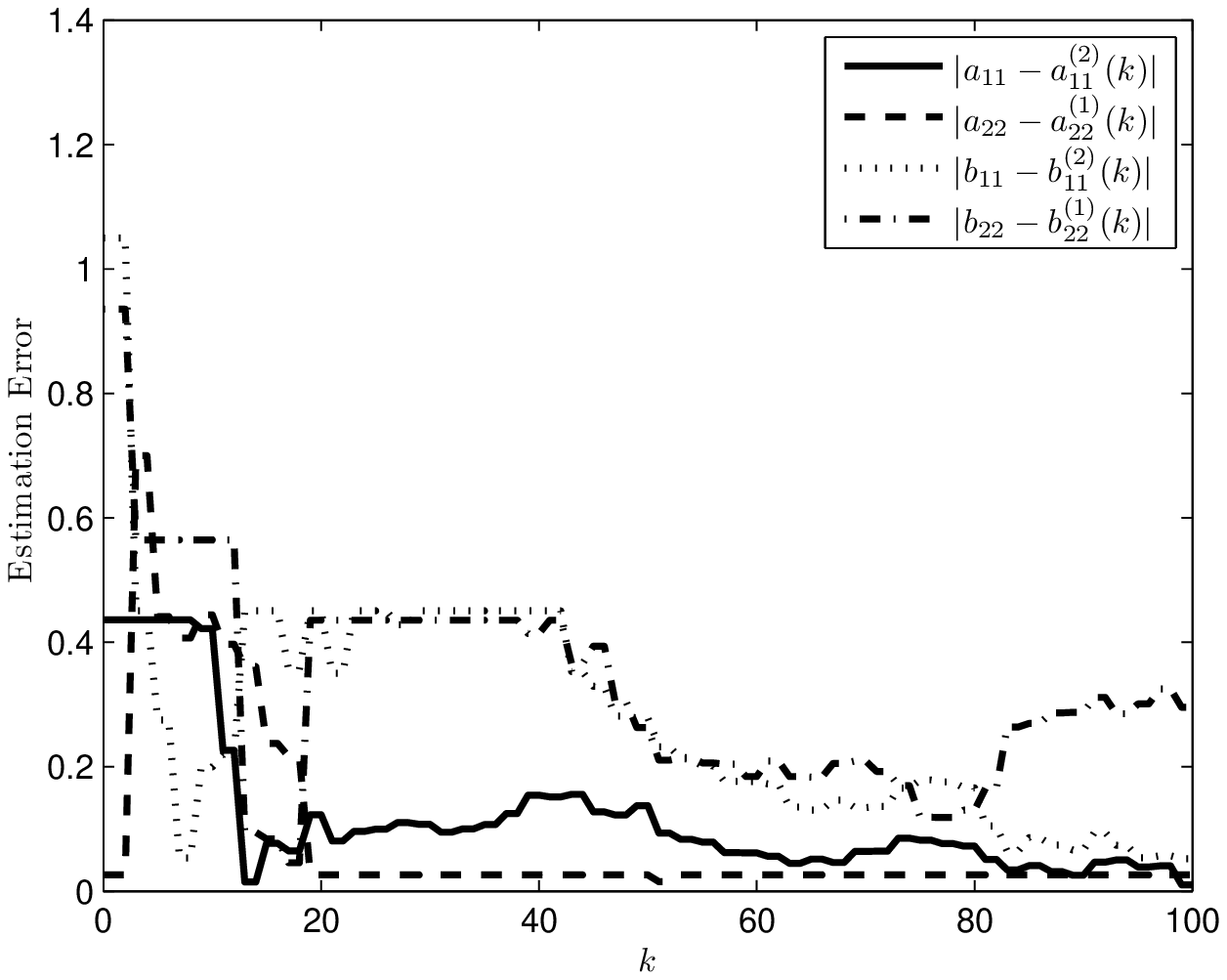}{figure2}{Estimation error of model parameters for the modified Campi-Kumar adaptive controller $\Gamma^*(P)$. }

\section{Conclusion} \label{sec:conclusion}
In this paper, as a generalization of earlier results in optimal control design with limited model information, we searched over the set of control design strategies that construct adaptive controllers. We found a minimizer of the competitive ratio both in average and supremum senses. We used the Campi--Kumar adaptive algorithm to setup an adaptive control design strategy that achieves a competitive ratio equal to one contrary to control design strategies that construct linear time-invariant control laws. This adaptive controller asymptotically achieves closed-loop performance equal to the optimal centralized controller with full model information. We illustrated the applicability of this adaptive controller on a vehicle platooning problem. As a future work, we suggest studying decentralized adaptive controllers.

\bibliographystyle{ieeetr}
\bibliography{ref}

\appendix

\section{Proof of Lemma~\ref{lem:aveleqsup}} \label{proof:lem:aveleqsup}
Let us assume, without loss of generality, that $r_\p^{\mathrm{sup}}(\Gamma)<\infty$ since otherwise, the desired inequality is trivially satisfied. First, note that using Theorem~2.10.1 in~\cite{friedman1970foundations}, function $J_P (\Gamma(P))/J_P (\mathbf{K}^*(P))$ is integrable on $\p$ since we assumed $r_\p^{\mathrm{sup}}(\Gamma)=\esssup J_P (\Gamma(P))/J_P (\mathbf{K}^*(P))<\infty$ (and $\p$ is a compact set due to Assumption~\ref{asm:1}). Then, using Theorem~2.7.1 in~\cite{friedman1970foundations}, we get
\begin{equation*}
\begin{split}
r_\p^{\mathrm{ave}}(\Gamma)=\int_{\xi\in\p} \frac{J_{\xi} (\Gamma({\xi}))}{J_{\xi} (\mathbf{K}^*({\xi}))} f({\xi})\; \mathrm{d} \xi 
&\leq \int_{ \xi \in\p} r_\p^{\mathrm{sup}}(\Gamma) f({\xi})\; \mathrm{d} \xi =r_\p^{\mathrm{sup}}(\Gamma).
\end{split}
\end{equation*}
This completes the proof.

\section{Proof of Lemma~\ref{tho:1}} \label{proof:tho:1}
Equations~(\ref{eqn:tho:1})--(\ref{eqn:tho:2.5}) are direct consequences of Theore-ms~2 and~3 in~\cite{campi:1890}. We start with proving~(\ref{eqn:tho:3}). To do so, let us prove $\|K(k)-\mathbf{L}(A ,B )\|> \rho$ implies that there exists at least an index $i$ such that $\|\mathbf{T}_iK^{(i)}(k)- \mathbf{T}_i \mathbf{L}(A ,B )\|> \rho/\sqrt{N}$. We can prove this fact by contradiction. Assume that there does not exists any index $i$ such that $\|\mathbf{T}_iK^{(i)}(k)-\mathbf{T}_i\mathbf{L}(A ,B )\|> \rho/\sqrt{N}$. Therefore, for all $1\leq i\leq N$, we have $\|\mathbf{T}_iK^{(i)}(k)-\mathbf{T}_i\mathbf{L}(A ,B )\|\leq \rho/\sqrt{N}$, and as a result, according to Theorem~1 in~\cite{springerlink:10.1007/BF01446925}, we get
$$
\|K(k)-\mathbf{L}(A ,B )\|^2\hspace{-.04in} \leq \hspace{-.04in}\sum_{i=1}^N \|\mathbf{T}_iK^{(i)}(k)-\mathbf{T}_i\mathbf{L}(A ,B )\|^2 \hspace{-.04in}\leq \hspace{-.04in}\rho^2\hspace{-.03in}.
$$
This is contradictory to the assumption that $\|K(k)-\mathbf{L}(A ,B )\|> \rho$. Hence, we proved the implication. Based on this property, it is easy to see that
\begin{equation} \label{eqn:proof:1}
\begin{split}
\sum_{t=0}^k &\chi(\|K(k)-\mathbf{L}(A ,B )\|> \rho) \leq \sum_{t=0}^k \sum_{i=1}^N \chi(\|\mathbf{T}_iK^{(i)}(k)-\mathbf{T}_i\mathbf{L}(A ,B )\|\hspace{-.05in}>\hspace{-.05in} \rho/\sqrt{N}).
\end{split}
\end{equation}
Now, note that $\|\mathbf{T}_iK^{(i)}(k)-\mathbf{T}_i\mathbf{L}(A ,B )\|> \rho/\sqrt{N}$ implies that $\|K^{(i)}(k)-\mathbf{L}(A ,B )\|> \rho/\sqrt{N}$. Thus, we get
\begin{equation} \label{eqn:proof:2}
\begin{split}
\sum_{t=0}^k \chi(\|\mathbf{T}_i&K^{(i)}(k)-\mathbf{T}_i\mathbf{L}(A ,B )\|> \rho/\sqrt{N}) \leq \sum_{t=0}^k \chi(\|K^{(i)}(k)-\mathbf{L}(A ,B )\|> \rho/\sqrt{N}).
\end{split}
\end{equation}
Substituting~(\ref{eqn:proof:2}) inside~(\ref{eqn:proof:1}), we get
\begin{equation*}
\begin{split}
\sum_{t=0}^k &\chi(\|K(k)-\mathbf{L}(A ,B )\|> \rho) \leq \sum_{t=0}^k \sum_{i=1}^N \chi(\|K^{(i)}(k)-\mathbf{L}(A ,B )\|> \rho/\sqrt{N}).
\end{split}
\end{equation*}
Now, using~(\ref{eqn:tho:2.5}), we can show that
$$
\sum_{t=0}^k \chi(\|K^{(i)}(k)-\mathbf{L}(A ,B )\|> \rho/\sqrt{N})\stackrel{as}{=}O(\mu(k)),
$$
for all $1\leq i\leq N$. Therefore, we have
$$
\sum_{t=0}^k \chi(\|K(k)-\mathbf{L}(A ,B )\|>\rho)\stackrel{as}{=}O(\mu(k)).
$$
Finally, note that the proof of~(\ref{eqn:tho:4}) is a direct result of applying~(\ref{eqn:tho:3}) to the proof of Theorem~5 in~\cite{campi:1890}. This concludes the proof.

\section{Proof of Lemma~\ref{lem:norm}} \label{proof:lem:norm}
First, note that
\begin{equation*}
\begin{split}
(X-Y)^\top P (X+Y)+(X+&Y)^\top P (X-Y)=2(X^\top P X-Y^\top P Y).
\end{split}
\end{equation*}
Hence, we get
\begin{equation*}
\begin{split}
2\|(X^\top P X-Y^\top P Y)\|&=\|(X-Y)^\top P (X+Y)+(X+Y)^\top P (X-Y)\|
\\&\leq \|(X-Y)^\top P (X+Y)\|+\|(X+Y)^\top P (X-Y)\|\\ & \leq 2\|P\|\|X-Y\|\|X+Y\|
\\ & \leq 2\|P\|\|X-Y\|(\|X\|+\|Y\|).
\end{split}
\end{equation*}
This concludes the proof.

\section{Proof of Corollary~\ref{cor:1}} \label{proof:cor:1}
First, notice that Theorem~4 implies $J_P(\Gamma^*(P))\stackrel{as}{=} J_P(\mathbf{K}^*(P))$, or equivalently $J_P(\Gamma^*(P))/J_P(\mathbf{K}^*(P))\stackrel{as}{=} 1$, for all $P\in\mathcal{P}\setminus\mathcal{N}$ (with $\mathcal{N}$ being a zero-measure set in the space $\bar{\mathcal{A}}\times\bar{\mathcal{B}}$). Therefore, by the definition of the essentially supremum operator (presented in Subsection~\ref{subsec:notation}), we get
$$
r_{\mathcal{P}}^{\mathrm{sup}} (\Gamma^*)=\esssup_{P\in\mathcal{P}} \frac{J_P (\Gamma^*(P))}{J_P (\mathbf{K}^*(P))}\stackrel{as}{=} 1.
$$
Now, applying Lemma~1 results in $1\leq r_{\mathcal{P}}^{\mathrm{ave}}\leq r_{\mathcal{P}}^{\mathrm{sup}}\stackrel{as}{=} 1$ and, hence, we have $r_{\mathcal{P}}^{\mathrm{ave}}\stackrel{as}{=} 1$.

\end{document}